\documentclass[12pt, twoside]{article}
\usepackage[tbtags]{amsmath}
\usepackage{amsfonts,amssymb, amsthm}
\usepackage{mathtext}
\usepackage{indentfirst}

\usepackage[T2A]{fontenc}
\usepackage[cp1251]{inputenc}
\usepackage[english,russian]{babel}

\newtheorem{lemma}{Лемма}
\newtheorem{theorem}{Теорема}
\begin{document}

\pagestyle{myheadings}
\markboth{Ф.~Гётце, Д.\,В.~Коледа, М.\,А.~Королёв}{О числе
квадратичных многочленов}

\title{О числе квадратичных многочленов с ограниченными дискриминантами}

\author{
Ф.~Гётце
\and
Д.\,В.~Коледа
\and
М.\,А.~Королёв
}

\maketitle

\renewcommand{\thefootnote}{}

\footnote{\emph{УДК} 511.42}

\footnote{\emph{Ключевые слова}: дискриминант многочлена, квадратичный многочлен.}

\renewcommand{\thefootnote}{\arabic{footnote}}
\setcounter{footnote}{0}

\begin{abstract}
В работе получена асимптотика для числа квадратичных
цело\-чис\-лен\-ных многочленов с ограниченной высотой и
дис\-кри\-ми\-нан\-та\-ми, не превосходящими заданной границы.

Библиография: 11 названий.
\end{abstract}

\section{Введение}

Пусть $p(x)=a_{n}x^{n}+\dots+a_{1}x+a_{0} \in \mathbb{Z}[x]$
многочлен степени $n$, $\displaystyle H(p) = \max_{0\le
i\le n} |a_{i}|$ -- высота многочлена $p$, и пусть
$\alpha_{1},\alpha_{2},\dots,\alpha_{n} \in \mathbb{C}$ -- корни
$p(x)$. Число
\begin{equation}\label{eq_01_Discrim}
 D(p)=a_{n}^{2n-2}\prod_{1\le i<j\le n}(\alpha_{i}-\alpha_{j})^{2}
\end{equation}
называется дискриминантом многочлена $p$. Дискриминант является
целочисленным многочленом от $n+1$ переменных коэффициентов
многочлена $p$ (см. \cite[гл.5, \S 33]{1}). Свойства дискриминанта
имеют многочисленные приложения в теории чисел и ал\-геб\-ре.

Так, например, поведение $D(p)$ тесно связано с задачей разделения
сопряженных алгебраических чисел \cite{2}-\cite{5}. Недавно в
\cite{6},\cite{7} были получены нижние оценки ко\-ли\-чест\-ва
многочленов с заданными высотами и ограниченными дискриминантами, а
также дискриминантами, делящимися на большую степень простого числа.
За\-ви\-си\-мость дискриминанта от коэффициентов многочлена имеет
сложный при $n=3$ и очень сложный при $n\!\ge\!4$ вид. Поэтому в
указанных выше работах использовалось определение
(\ref{eq_01_Discrim}).

С другой стороны, при малых $n$ средствами аналитической теории
чисел удается получать более сильные результаты. Так в \cite{8} была
найдена верхняя оценка для количества целочисленных многочленов
третьей степени с заданной границей для абсолютной величины
дискриминанта.

Пусть $n\ge 2$ -- целое число, $Q>1$, $v\ge 0$. Определим множества
\begin{equation*}
\mathcal{P}_{n}(Q)=\{p(x)\in \mathbb{Z}[x] : \deg p = n, \ H(p) \le
Q\},
\end{equation*}
\begin{equation*}
\mathcal{P}_{n}(Q,v)=\{p(x)\in \mathcal{P}_{n}(Q) : |D(p)|\le
\gamma_n Q^{2n-2-2v}\},
\end{equation*}
где
$$
\gamma_n = \sup_{\substack{
p \in \mathbb{Z}[x]\\
\deg p = n }} \frac{|D(p)|}{(H(p))^{2n-2}}.
$$
Заметим, что условие $H(p)\le Q$ влечет неравенство $|D(p)|\le
\gamma_n Q^{2n-2}$. В частности, при $n = 2$ и $p(x) = ax^2 + bx +
c$ имеем: $D(p) = b^{2}-4ac$, откуда $|D(p)|\le 5Q^{2}$, причем
равенство достигается на многочлене $Q(x^{2}+x-1)$. Таким образом,
$\gamma_{2} = 5$.

В \cite{6} и \cite{9} было доказано, что
\begin{equation}\label{Lab_02}
\# \mathcal{P}_n(Q, v) \gg Q^{n+1-2v}
\end{equation}
при $0 < v < \frac{1}{2}$. Из эвристических соображений можно
предположить, что с точ\-ностью до постоянного множителя правая
часть соотношения (\ref{Lab_02}) будет также и верхней границей для
величины $\# \mathcal{P}_n(Q, v)$. Однако в \cite{8} было доказано,
что $\# \mathcal{P}_3(Q,v) \ll Q^{\,4-5v/3}$. Незначительное
видоизменение рассуждений из работы \cite{8} приводит к соотношению
$\# \mathcal{P}_3(Q,v) \asymp Q^{\,4-5v/3}$. Это показывает сложную
за\-ви\-си\-мость различных характеристик (сумм, интегралов, мер и
т.д.) многочлена от его ко\-эф\-фи\-ци\-ен\-тов. Пожалуй, впервые
данное обстоятельство было замечено в работе \cite{10} при решении
проблемы Терри.

Возникает естественный вопрос: когда оценка (\ref{Lab_02}) является
правильной по по\-ряд\-ку? В настоящей работе на этот вопрос дается
утвердительный ответ в случае, когда $n = 2$, а $v$ - фиксированное
число с условием $0<v<\frac{3}{4}$. Именно, пусть
$N(Q,D)$ -- количество целочисленных многочленов второй степени,
высота которых не превосходит $Q$, а дискриминант ограничен по
абсолютной величине числом $D$ (так что $\#\mathcal{P}_{2}(Q,v) =
N(Q,5Q^{2-2v})$). Тогда имеет место

\begin{theorem}\label{Theorem_1}
При любых $D$ и $Q$ с условием $1\le D\le Q^{2}/2$ справедливо
равенство
\[
N(Q,D)\,=\,\kappa QD\,+\,O\left(D^{3/2}\ln Q + (Q \ln Q)^{3/2}\right),
\]
где $\kappa = 4(\ln{2}+1) =  6.772588\ldots$.
\end{theorem}

\textsc{Замечание.} Несложно проверить, что формула для $N(Q,D)$
будет асимп\-то\-ти\-чес\-кой, если
\[
Q^{1/2} (\ln{Q})^{3/2}\,\ll\,D\,\ll\,\biggl(\frac{Q}{\ln{Q}}\biggr)^{2}.
\]
Эти условия выполнены, если $D = Q^{2-2v}$, где $v$ -- фиксированное
число из про\-ме\-жут\-ка $\bigl(0,\frac{3}{4}\bigr)$.

В статье используются стандартные теоретико-числовые обозначения:
$[x]$ и $\{x\} = x-[x]$ -- соответственно, целая и дробная части
числа $x$, $\lceil x \rceil$ -- ближайшее к $x$ целое число, не меньшее чем $x$, $\|x\| = \min{(\{x\},1-\{x\})}$ -- расстояние от $x$ до
ближайшего целого числа; $\tau(n)$ -- число делителей $n$ или, что
то же, число решений уравнения $xy = n$ в натуральных $x$ и $y$;
$\theta, \theta_{1},\ldots$ -- комплексные числа, по модулю не
превосходящие единицы и в разных соотношениях, вообще говоря,
разные; $(a,b)=\nod(a,b)$ --- наибольший общий делитель целых чисел $a$ и $b$.

\section{Вспомогательные утверждения}

В данном разделе собраны вспомогательные утверждения, необходимые для доказательства основной теоремы.

\begin{lemma}\label{lmm1}
Пусть
\[
\alpha = \frac{a}{q} + \frac{\theta}{q^2}, \qquad (a,q) = 1, \qquad q\ge 1, \qquad |\theta|\le 1.
\]
Тогда при любом $\beta$, $U>0$, $P\ge 1$ имеем
\[
\sum_{x=1}^P \min\left(U, \frac{1}{\|\alpha x + \beta\|}\right) \le 6 \left(\frac{P}{q} + 1\right) (U + q \ln q).
\]
\end{lemma}

Доказательство см. стр. 94--95 в \cite[\S 2 гл. VI]{13}.

\begin{lemma}\label{lmm2}
Пусть $\nod(a,m)=1$, $1\le N\le m$. Тогда для любых $m$, больших некоторой абсолютной постоянной $m_0$, справедливо неравенство
\[
\left|\sum_{x=1}^N {e^{2\pi i \frac{ax^2}{m}}}\right| \le 5 \sqrt{m \ln m}\,.
\]
\end{lemma}

\textsc{Доказательство.}
Пусть $S = \sum_{x=1}^N e^{2\pi i \frac{a x^2}{m}}$. Тогда
\begin{multline*}
|S|^2 = \sum_{x,y=1}^N e^{2\pi i \frac{a\left(y^2-x^2\right)}{m}} = \sum_{1\le x\le N} \sum_{1\le x+h\le N} e^{2\pi i \frac{a\left((x+h)^2-x^2\right)}{m}} =\\
= \sum_{1\le x\le N} \sum_{1-x\le h\le N-x} e^{2\pi i \frac{2axh + ah^2}{m}} = \sum_{|h|<N} e^{2\pi i \frac{ah^2}{m}} \sum_{X_1\le x\le X_2} e^{2\pi i \frac{2ahx}{m}},
\end{multline*}
где $X_1 = \max(1,1-h)$, $X_2 = \min(N,N-h)$, $X_2 - X_1 < N$.

Выделяя слагаемые, отвечающие $h=0$, получим:
\begin{equation*}
|S|^2 \le N + \sum_{0<|h|<N} S_h, \qquad
S_h = \left|\sum_{X_1\le x\le X_2} e^{2\pi i \frac{2ahx}{m}}\right|.
\end{equation*}
Имеем, далее
\[
S_h\le \min\left(X_2-X_1+1, \frac{1}{\left|\sin \frac{2\pi a h}{m}\right|}\right) \le \min\left(N, \frac{1}{2\left\|\frac{2ah}{m}\right\|}\right).
\]
Замечая, что суммы $S_h$ и $S_{-h}$ оцениваются одинаково, найдём:
\[
\sum_{0<|h|<N} S_h \le 2 \sum_{1\le h < N} \min\left(N, \frac{1}{2\left\|\frac{2ah}{m}\right\|}\right) \le 2 \sum_{1\le h < 2N} \min\left(N, \frac{1}{2\left\|\frac{ah}{m}\right\|}\right).
\]
Полагая $\alpha = \frac{a}{m}$, $\beta = 0$, $P=U=2N$ в лемме \ref{lmm1} и замечая, что $N\le m$, при $m\ge m_0$ получим:
\begin{multline*}
\sum_{0<|h|<N} S_h \le \sum_{1\le h< 2N} \min\left(2N, \left\| \frac{ah}{m}\right\|^{-1}\right) \le \\
\le  6 \left(\frac{2N}{m}+1\right) (2N+m \ln m) < 20 m \ln m.
\end{multline*}
Таким образом, получаем
\[
|S|^2 \le N + 20 m \ln m \le 25 m \ln m.
\]
Лемма доказана.

\begin{lemma}\label{lmm3}
Пусть $m\in\mathbb{N}$, $A_1$, $A_2$, $B_1$, $B_2$ --- целые числа, $A_1 \le A_2$, $B_1\le B_2$.
Тогда
\[
\sum_{A_1\le a\le A_2} \sum_{\substack{B_1 \le q\le B_2\\ q^2\equiv a\!\!\!\!\!\pmod{m}}} \!\!1
\le
\left\lceil \frac{A_2 - A_1 + 1}{m} \right\rceil (B_2-B_1+1).
\]
\end{lemma}

\textsc{Доказательство.}
Заметим, что
\[
\sum_{0 \le a < m} \sum_{\substack{B_1 \le q\le B_2\\ q^2\equiv a\!\!\!\!\!\pmod{m}}} 1 =
B_2 - B_1 + 1.
\]
Каждый квадратичный вычет по модулю $m$ встречается в промежутке $A_1\le a\le A_2$ не более $\left\lceil \frac{A_2 - A_1 + 1}{m} \right\rceil$ и не менее
$\left[\frac{A_2 - A_1 + 1}{m} \right]$ раз. Лемма доказана.

\section{Доказательство теоремы}

Подсчитаем число троек целых чисел $(q,n,r)$ с условиями $-Q\le q,n,r \le Q$ и $-D\le q^2-4nr \le D$ (предполагаем, что $1\le D\le Q^2/2$).

Если $q=0$, то последнее ограничение на $q,n,r$ принимает вид $-\frac{D}{4}\le nr \le \frac{D}{4}$. Если одно из чисел $n,r$ равно нулю, то для второго имеем не более $2Q+1$ возможностей. Если $nr\ne 0$, то число таких пар $(n,r)$ не превосходит
\[
4 \sum_{1\le k\le D/4} \tau(k) \ll D\ln D \ll D\ln Q.
\]
Итак, количество троек $(q,n,r)$, где $q=0$, не превосходит по порядку $Q+D\ln Q$.

Рассмотрим теперь те тройки $(q,n,r)$, где $q\ne 0$, но $nr = 0$. Тогда $1\le q^2\le D$, откуда $1\le |q|\le \sqrt{D}$. Так как для ненулевого числа из пары $(n,r)$ имеется не более $2Q+1$ возможностей, то число троек с условием $q\ne 0$, $nr=0$, не превосходит по порядку $Q\sqrt{D}$.

Оставшиеся тройки отвечают условию $qnr\ne 0$. Очевидно, число троек с положительным $q$ совпадает с числом троек с отрицательным $q$. Далее, число троек с $n<0$, $r<0$ совпадает с числом троек с $n>0$, $r>0$, а число троек с $n<0$, $r>0$ равно числу троек с $n>0$, $r<0$. Следовательно,
\[
N(Q,D) = 4(N_1+N_2) + O(Q\sqrt{D}+D\ln Q),
\]
где $N_1$ --- число троек с условиями $1\le q,n,r \le Q$, $-D\le q^2-4nr\le D$,
а $N_2$ --- число троек с условиями $1\le q,n,r \le Q$, $-D\le q^2 + 4nr \le D$
(последнее неравенство, очевидно, можно заменить следующим: $1\le q^2+4nr\le D$).

Ясно, что $N_2$ не превосходит числа троек с условиями $1\le q^2 \le D$, $1\le 4nr\le D$, $1\le n,r \le Q$  или, что то же, $1\le q\le \sqrt{D}$, $1\le nr\le D/4$, $1\le n,r\le Q$, и, таким образом, оценивается сверху величиной порядка $D^{3/2} \ln D\ll D^{3/2} \ln Q$. Итак,
\[
N(Q,D) = 4N_1 + O\!\left(Q\sqrt{D}+D^{3/2}\ln Q\right).
\]

Вычислим $N_1$. Зафиксируем целое $a$, $|a|\le D$, и обозначим через $N_1(a)$ число тех троек $(q,n,r)$, для которых $q^2-4nr=a$ при условиях $1\le q,n,r \le Q$.
(При $a\not\equiv 0\pmod{4}$ и $a\not\equiv 1\pmod{4}$ заведомо имеем $N_1(a)=0$.
Действительно, необходимо $q^2\equiv a\pmod{4}$. Но $q^2\equiv 0\pmod{4}$ или $q^2\equiv 1\pmod{4}$.)
Далее имеем:
\begin{equation*}
N_1(a) = \sum_{1\le q\le Q} \sum_{\substack{1\le n,r \le Q\\q^2-a=4nr}} 1 =
\sum_{\substack{1\le n,r\le Q\\nr\le \frac{Q^2-a}{4}}} \sum_{\substack{1\le q\le Q\\q^2-a=4nr}} 1 =
2 S_1(a) - S_2(a),
\end{equation*}
где
\begin{gather*}
S_1(a) := \sum_{1\le n \le \frac{\sqrt{Q^2-a}}{2}} \sum_{\substack{1\le r\le Q \\r\le \frac{Q^2-a}{4n}}} \sum_{\substack{1\le q\le Q\\q^2-a=4nr}} 1,\\
S_2(a) := \sum_{1\le n,r\le \frac{\sqrt{Q^2-a}}{2}} \sum_{\substack{1\le q\le Q\\q^2-a=4nr}} 1.
\end{gather*}

Предположим, что для некоторых $n$ и $q$, $1\le n\le \frac{\sqrt{Q^2-a}}{2}$,
$1\le q\le Q$, выполняется сравнение $q^2\equiv a\pmod{4n}$. Положив $r=\frac{q^2-a}{4n}$, получим $r\le \frac{Q^2-a}{4n}$. Чтобы $r$ удовлетворяло условию $1\le r\le Q$,
необходимо потребовать: $1\le \frac{q^2-a}{4n} \le Q$, т.е. $4n+a\le q^2\le 4nQ+a$.
Аналогично, для того, чтобы $r$ удовлетворяло условию $1\le r\le \frac{\sqrt{Q^2-a}}{2}$,
необходимо выполнение неравенства $1\le \frac{q^2-a}{4n}\le\frac{\sqrt{Q^2-a}}{2}$,
т.е. $4n+a\le q^2\le 2n\sqrt{Q^2-a}+a$.
Так получим
\[
\begin{gathered}
S_1(a)\  = \sum_{1\le n\le\frac{\sqrt{Q^2-a}}{2}} \sum_{\substack{1\le q\le Q\\ 4n+a\le q^2 \le 4nQ+a\\ q^2\equiv a\!\!\!\!\!\pmod{4n}}} 1\,,\\
S_2(a) \ = \sum_{1\le n\le\frac{\sqrt{Q^2-a}}{2}} \sum_{\substack{1\le q\le Q\\ 4n+a\le q^2 \le 2n\sqrt{Q^2-a}+a\\ q^2\equiv a\!\!\!\!\!\pmod{4n}}} \!\!\!\!1\,.
\end{gathered}
\]

Таким образом, получаем, что $N_1 = \sum_{|a|\le D} N_1(a) = 2S_1-S_2$, где
\[
S_1 = \sum_{|a|\le D} S_1(a), \qquad S_2 = \sum_{|a|\le D} S_2(a).
\]

Далее преобразуем суммы $S_1$ и $S_2$ с тем, чтобы границы изменения $n$ и $q$ не зависели от $a$. А именно, обозначим через $T_1$ вклад в $S_1$ слагаемых, отвечающих $1\le a\le D$, $\frac{\sqrt{Q^2-a}}{2} < n \le \frac{Q}{2}$, а через $T_2$ --- вклад слагаемых, соответствующих $-D\le a\le -1$, $\frac{Q}{2} < n\le \frac{\sqrt{Q^2-a}}{2}$.
Так что $S_1 = S_1' - T_1 + T_2$, где
\begin{equation}
S_1' \ = \sum_{|a|\le D} \sum_{1\le n\le \frac{Q}{2}} \sum_{\substack{1\le q\le Q\\ 4n+a\le q^2\le 4nQ+a\\ q^2\equiv a\!\!\!\!\!\pmod{4n}}} 1\,.
\end{equation}

Далее при оценках внутренних сумм будем пользоваться леммой \ref{lmm3}.

Пусть $|a|\le D\le Q^2/2$. Применяя теорему о среднем (теорему Лагранжа), получаем
\[
\frac{\sqrt{Q^2-a}}{2} = \frac{Q}{2} \sqrt{1-\frac{a}{Q^2}}=\frac{Q}{2}-\frac{1}{4\sqrt{1-\xi}}\frac{a}{Q},
\]
где $\xi$ --- некоторая точка между нулём и значением $a/Q^2$.

В силу условия $|a|\le D\le Q^2/2$ имеем $|\xi|\le 1/2$, и, следовательно, между $\frac{\sqrt{Q^2-a}}{2}$ и $Q/2$
лежит не более $\frac{D}{2\sqrt{2} Q} + 1\ll D/Q+1$ целых чисел. Итак, имеем
\begin{multline*}
T_1 := \sum_{1\le a\le D} \sum_{\frac{\sqrt{Q^2-a}}{2}< n \le \frac{Q}{2}}
\sum_{\substack{1\le q\le Q\\4n+a\le q^2\le 4nQ+a\\q^2\equiv a\!\!\!\!\pmod{4n}}} \!\!\!1 \ \le
 \sum_{\frac{\sqrt{Q^2-D}}{2}< n \le \frac{Q}{2}} \sum_{1\le a\le D}
\sum_{\substack{1\le q\le Q\\4n+1\le q^2\le 4nQ+D\\q^2\equiv a\!\!\!\!\pmod{4n}}} \!\!\!1 \ \ll \\
\ll \sum_{\frac{\sqrt{Q^2-D}}{2}< n \le \frac{Q}{2}} \left(\frac{D}{n}+1\right) Q \ \ll \
Q \left(D/Q+1\right)^2 \ \ll \ D^2/Q + Q.
\end{multline*}

Аналогично оценивается и сумма $T_2$,
\[
T_2 := \sum_{-D\le a\le -1} \sum_{\frac{Q}{2}< n \le\frac{\sqrt{Q^2-a}}{2}}
\sum_{\substack{1\le q\le Q\\4n+a\le q^2\le 4nQ+a\\q^2\equiv a\!\!\!\!\pmod{4n}}} \!\!\!1 \ \ll  D^2/Q + Q.
\]

Определим для суммы $S_2$ величины $T_3$ и $T_4$ подобно $T_1$ и $T_2$. Так получим:
\begin{multline*}
T_3 := \sum_{1\le a\le D} \sum_{\frac{\sqrt{Q^2-a}}{2}< n \le \frac{Q}{2}}
\sum_{\substack{1\le q\le Q\\4n+a\le q^2\le 2n\sqrt{Q^2-a}+a\\q^2\equiv a\!\!\!\!\pmod{4n}}} \!\!\!1 \ \le\\
\le
\sum_{\frac{\sqrt{Q^2-D}}{2}< n \le \frac{Q}{2}} \sum_{1\le a\le D}
\sum_{\substack{1\le q\le Q\\4n+1\le q^2\le 4n\frac{Q}{2}+D\\q^2\equiv a\!\!\!\!\pmod{4n}}} \!\!\!1 \ \ \ll \ \ 
D^2/Q + Q,
\end{multline*}
\[
T_4 := \sum_{-D\le a\le -1} \sum_{\frac{Q}{2}< n \le\frac{\sqrt{Q^2-a}}{2}}
\sum_{\substack{1\le q\le Q\\4n+a\le q^2\le 2n\sqrt{Q^2-a}+a\\q^2\equiv a\!\!\!\!\pmod{4n}}} \!\!\!1\ \ll \
D^2/Q + Q.
\]

Итак, $N_1 = 2 S_1' - S_2' + O(D^2/Q + Q)$, где
\[
S_1' := S_1 + T_1 - T_2 = \sum_{|a|\le D} \sum_{1\le n\le \frac{Q}{2}} \sum_{\substack{1\le q\le Q\\4n+a\le q^2 \le 4nQ+a\\q^2\equiv a\!\!\!\!\!\pmod{4n}}} 1,
\]
\[
S_2' := S_2 + T_3 - T_4 = \sum_{|a|\le D} \sum_{1\le n\le \frac{Q}{2}} \sum_{\substack{1\le q\le Q\\4n+a\le q^2 \le 2n\sqrt{Q^2-a}+a\\q^2\equiv a\!\!\!\!\!\pmod{4n}}} 1.
\]

Покажем теперь, что в $S_1'$ и $S_2'$ границы изменения $q$ можно сделать не зависящими от $a$. Пусть, скажем,
$T_5$ --- вклад слагаемых $S_1'$, отвечающих $1\le a\le D$ и $4n \le q^2 < 4n + a$, $T_6$ --- вклад слагаемых,
соотвествующих $-D\le a\le-1$,
$4n+a\le q^2 < 4n$ (в каждом из случаев $1\le n\le Q/2$). Тогда
\begin{multline*}
T_5 := \sum_{1\le a\le D} \sum_{1\le n\le \frac{Q}{2}} \sum_{\substack{1\le q\le Q\\4n\le q^2 < 4n+a\\q^2\equiv a\!\!\!\!\!\pmod{4n}}}  1 \  \le \ 
\sum_{1\le n\le \frac{Q}{2}} \sum_{1\le a\le D} \sum_{\substack{4n\le q^2\le 4n+D\\q^2\equiv a\!\!\!\!\!\pmod{4n}}} 1 \ll \\
\ll
\sum_{1\le n\le\frac{Q}{2}} \left(\frac{D}{n} + 1\right)\left(\frac{D}{\sqrt{n+D}}+1\right)
\ll
\sum_{1\le n\le \frac{Q}{2}} \left(\frac{D^{3/2}}{n} + \frac{D}{\sqrt{n}} + 1\right) \ll \\
\ll
D^{3/2} \ln Q + D\sqrt{Q} + Q.
\end{multline*}

Аналогично получается и оценка
\begin{multline*}
T_6 := \sum_{-D\le a\le -1} \sum_{1\le n\le \frac{Q}{2}} \sum_{\substack{1\le q\le Q\\4n+a\le q^2 < 4n\\q^2\equiv a\!\!\!\!\!\pmod{4n}}} 1 \ \ll \\
\ll \sum_{1\le n\le\frac{Q}{2}} \left(\frac{D}{n}+1\right)\left(\sqrt{4n}-\sqrt{\max(4n-D, 0)}+1\right)
\ \ll \ 
D^{3/2} + D\sqrt{Q} + Q.
\end{multline*}

Пусть, далее, $T_7$ --- вклад в $S_1'$ слагаемых, отвечающих $1\le a\le D$, $4nQ<q^2\le 4nQ+a$,
$T_8$ --- от слагаемых, отвечающих $-D\le a\le -1$, $4nQ+a < q^2 \le 4nQ$
(везде $1\le n \le Q/2$). Тогда
\begin{multline*}
T_7 := \sum_{1\le a\le D}\sum_{1\le n\le \frac{Q}{2}} \sum_{\substack{1\le q\le Q\\4nQ<q^2\le 4nQ+a\\q^2\equiv a\!\!\!\!\!\pmod{4n}}} 1 \ \le\\
\le \sum_{1\le n\le \frac{Q}{2}} \sum_{1\le a\le D} \sum_{\substack{4nQ<q^2\le 4nQ+D\\q^2\equiv a\!\!\!\!\!\pmod{4n}}} 1 \ll
\sum_{1\le n\le Q/2} \left(\frac{D}{n}+1\right) \left(\frac{D}{\sqrt{nQ+D}}+1\right) \ll \\
\ll
\sum_{1\le n\le Q/2} \left(\frac{D^{3/2}}{n}+\frac{D}{n}+\frac{D}{\sqrt{nQ}}+1\right)\ \ll \ 
D^{3/2}\ln Q + Q,
\end{multline*}
и, аналогично,
\begin{multline*}
T_8 := \sum_{-D\le a\le -1} \sum_{1\le n\le\frac{Q}{2}} \sum_{\substack{1\le q \le Q\\ 4nQ+a < q^2 \le 4nQ\\q^2\equiv a\!\!\!\!\!\pmod{4n}}} 1 \ \ll \\
\ll \sum_{1\le n\le \frac{Q}{2}} \left(\frac{D}{n}+1\right) \left(\sqrt{4nQ}-\sqrt{\max(4nQ-D, 0)} + 1\right) \ 
\ll \ D^{3/2} + D\ln Q + Q.
\end{multline*}

Определим для суммы $S_2'$ величины $T_9,\ldots,T_{12}$ по аналогии с $T_5,\ldots,T_8$.
Замечая, что $T_9 = T_5$, $T_{10} = T_6$,
будем иметь
\[
T_9 \ll D^{3/2}\ln Q + D\sqrt{Q}+Q, \qquad
T_{10} \ll D^{3/2} + D\sqrt{Q}+Q,
\]
\begin{multline*}
T_{11} := \sum_{1\le a\le D} \sum_{1\le n\le \frac{Q}{2}} \sum_{\substack{2nQ< q^2 \le 2n\sqrt{Q^2-a}+a\\ q^2\equiv a\!\!\!\!\!\pmod{4n}}}  1 \ \le \
\sum_{1\le n\le \frac{Q}{2}} \sum_{1\le a\le D} \sum_{\substack{2nQ< q^2 \le 2nQ+D\\ q^2\equiv a\!\!\!\!\!\pmod{4n}}}  1 \ \ll \\
\ll
\sum_{1\le n\le \frac{Q}{2}} \left(\frac{D}{n}+1\right) \left(\frac{D}{\sqrt{2nQ+D}}+1\right) \  \ll \ 
D^{3/2} \ln Q + Q,
\end{multline*}
\begin{multline*}
T_{12} := \sum_{-D\le a\le -1} \sum_{1\le n\le \frac{Q}{2}} \sum_{\substack{1\le q\le Q\\2n\sqrt{Q^2-a}+a < q^2 \le 2nQ\\ q^2\equiv a\!\!\!\!\!\pmod{4n}}} 1 \ \le
\sum_{-D\le a\le -1} \sum_{1\le n\le \frac{Q}{2}} \sum_{\substack{1\le q\le Q\\2nQ-D\le q^2 \le 2nQ\\ q^2\equiv a\!\!\!\!\!\pmod{4n}}} 1 \ \ll \\
\ll
\sum_{1\le n\le \frac{Q}{2}} \left(\frac{D}{n}+1\right)\left(\sqrt{2nQ}-\sqrt{\max(2nQ-D,0)}+1\right)
\ll D^{3/2} + D\ln Q + Q.
\end{multline*}

Окончательно находим
\[
N_1 = 2S_1'' - S_2'' + O\left(Q+D\sqrt{Q} + D^2/Q + D^{3/2} \ln Q \right),
\]
\[
\begin{gathered}
S_1'' := S_1' + T_5-T_6-T_7+T_8 = \sum_{1\le n\le \frac{Q}{2}}\ \sum_{|a|\le D} \sum_{\substack{1\le q\le Q\\4n\le q^2\le 4nQ\\q^2\equiv a\!\!\!\!\!\pmod{4n}}} 1\ , \\
S_2'' := S_2' + T_9-T_{10}-T_{11}+T_{12} = \sum_{1\le n\le \frac{Q}{2}}\ \sum_{|a|\le D} \sum_{\substack{1\le q\le Q\\4n\le q^2\le 2nQ\\q^2\equiv a\!\!\!\!\!\pmod{4n}}} 1\ .
\end{gathered}
\]

Вычислим $S_1''$. Имеем
\[
S_1'' := \sum_{1\le n\le \frac{Q}{2}}\ \sum_{|a|\le D} \sum_{\substack{2\sqrt{n}\le q\le 2\sqrt{nQ}\\1\le q\le Q\\q^2\equiv a\!\!\!\!\!\pmod{4n}}} 1\ .
\]

Если $1\le n\le Q/4$, то верхней границей $q$ во внутренней сумме будет величина $2\sqrt{nQ}$; в противном случае --- величина $Q$. Значит,
\[
S_1'' = V_1 + V_2,
\]
где
\[
\begin{gathered}
V_1 = \sum_{1\le n\le \frac{Q}{4}} \sum_{|a|\le D} \sum_{\substack{2\sqrt{n}\le q\le 2\sqrt{nQ}\\q^2\equiv a\!\!\!\!\!\pmod{4n}}} 1, \\
V_2 = \sum_{\frac{Q}{4}< n\le \frac{Q}{2}} \sum_{|a|\le D} \sum_{\substack{2\sqrt{n}\le q\le Q\\q^2\equiv a\!\!\!\!\!\pmod{4n}}} 1.
\end{gathered}
\]

Имеем, далее:
\begin{multline*}
V_1 =
\sum_{1\le n\le\frac{Q}{4}}\ \sum_{|a|\le D}\ \frac{1}{4n}\sum_{-2n<c\le 2n} \sum_{2\sqrt{n}\le q\le 2\sqrt{nQ}} e^{2\pi i\frac{c(q^2-a)}{4n}} = \\
= \sum_{1\le n\le\frac{Q}{4}}\ \frac{1}{4n}\sum_{-2n<c\le 2n} \left(\sum_{|a|\le D} e^{-2\pi i\frac{ac}{4n}}\right) \cdot
\left(\sum_{2\sqrt{n}\le q\le 2\sqrt{nQ}} e^{2\pi i\frac{cq^2}{4n}}\right) = \\
= \sum_{1\le n\le\frac{Q}{4}} \frac{1}{4n} \left(2[D]+1\right)\left(\left[2\sqrt{nQ}\right]-\left[2\sqrt{n}\right]+1\right) + R_1,
\end{multline*}
где первое слагаемое (сумма по $n$) отвечает вкладу от $c=0$, а слагаемое $R_1$ --- вкладу $c\ne 0$.

Первое слагаемое преобразуется к виду
\begin{multline*}
\frac14 \sum_{1\le n\le \frac{Q}{4}} \frac{1}{n} \left(2D+O(1)\right)\left(2\sqrt{nQ}+O(\sqrt{n})\right) 
= \\
= \frac14 \sum_{1\le n\le\frac{Q}{4}} \frac1n \left(4D\sqrt{nQ}+ O(D\sqrt{n})+ O(\sqrt{nQ})\right) = \\
= D\sqrt{Q} \sum_{1\le n\le\frac{Q}{4}} \frac{1}{\sqrt{n}} + O\left(D\sqrt{Q}+Q\right) = 
D\sqrt{Q} \left(2\sqrt{\frac{Q}{4}} + O(1)\right) + O\left(D\sqrt{Q}+Q\right) =\\
= DQ + O\left(D\sqrt{Q}+Q\right).
\end{multline*}

Пусть, далее, $\nod(a,m)=\delta \ge 1$, $a=\delta a_1$, $m=\delta m_1$, где $\nod(a_1,m_1)=1$. Тогда при любом $X>1$ имеем:
\begin{multline*}
\sum_{x=1}^X e^{2\pi i\frac{ax^2}{m}}=\sum_{x=1}^X e^{2\pi i\frac{a_1 x^2}{m_1}}=
\sum_{s=0}^{\left[\frac{X}{m_1}\right]-1} \sum_{x=sm_1+1}^{(s+1)m_1} e^{2\pi i\frac{a_1 x^2}{m_1}} +
\sum_{x=\left[\frac{X}{m_1}\right]m_1+1}^{X} e^{2\pi i\frac{a_1 x^2}{m_1}}=\\
=
\left[\frac{X}{m_1}\right] \sum_{x=1}^{m_1} e^{2\pi i \frac{a_1 x^2}{m_1}}\ +
\sum_{x=\left[\frac{X}{m_1}\right]m_1+1}^{X} e^{2\pi i\frac{a_1 x^2}{m_1}}
\end{multline*}

Пользуясь оценкой неполной суммы Гаусса из леммы \ref{lmm2}, получаем
\[
\left|\sum_{x=1}^X e^{2\pi i\frac{ax^2}{m}}\right| \ll \left(\frac{X}{\sqrt{m}}\sqrt{(a,m)} + \sqrt{\frac{m}{(a,m)}}\right) \sqrt{\ln m}.
\]

Замечая теперь, что при $0<|c| \le 2n$
\[
\left|\sum_{|a|\le D} e^{-2\pi i \frac{ac}{4n}}\right| = \left|\frac{\sin\left(\frac{2 \pi c ([D]+1/2)}{4n}\right)}{\sin\left(\frac{\pi c}{4n}\right)}\right| \ll \frac{n}{|c|}.
\]
Таким образом, имеем
\begin{multline*}
R_1 = \frac14 \sum_{1\le n\le\frac{Q}{4}} \frac1n \sum_{\substack{-2n < c \le 2n\\c\ne 0}}
\left(\sum_{|a|\le D} e^{-2\pi i \frac{ac}{4n}}\right) \cdot
\left(\sum_{2\sqrt{n}\le q\le 2\sqrt{nQ}} e^{2\pi i \frac{c q^2}{4n}}\right) \ll \\
\ll
\sum_{1\le n\le \frac{Q}{4}}\frac1n \sum_{1\le c\le 2n} \frac{n}{c} \left(\frac{\sqrt{nQ}}{\sqrt{n}}\sqrt{(c,4n)} + \sqrt{\frac{n}{(c,4n)}}\right)\sqrt{\ln n} \ll \\
\ll
\sum_{1\le n\le\frac{Q}{4}} \sum_{1\le c\le 2n} \frac1c \sqrt{Q \ln Q}\sqrt{(c,4n)} \ll
\sqrt{Q \ln Q} \sum_{1\le n\le Q} \sum_{c=1}^n \frac{\sqrt{(c,n)}}{c} \ll \\
\ll
\sqrt{Q \ln Q} \sum_{1\le n\le Q} \sum_{\delta|n} \sum_{\substack{1\le c\le n\\(c,n)=\delta}} \frac{\sqrt{\delta}}{c}.
\end{multline*}

Полагая в последней сумме $c=c_1\delta$, $n=n_1\delta$, где $(c_1,n_1)=1$, получим
\begin{multline*}
R_1 \ll \sqrt{Q \ln Q} \sum_{n\le Q} \sum_{\delta| n} \frac{1}{\sqrt{\delta}} \sum_{1\le c_1\le n_1} \frac{1}{c_1} \ll
\sqrt{Q} (\ln Q)^{3/2} \sum_{n\le Q} \sum_{\delta|n}\frac{1}{\sqrt{\delta}} \ll \\
\ll
\sqrt{Q} (\ln Q)^{3/2} \sum_{\delta\le Q} \frac{1}{\sqrt{\delta}} \sum_{1\le m\le \frac{Q}{\delta}} 1 \ll
\sqrt{Q} (\ln Q)^{3/2} \sum_{\delta\le Q} \frac{Q}{\delta \sqrt{\delta}} \ll (Q \ln Q)^{3/2}.
\end{multline*}

Итак, $V_1 = DQ + O\left(D\sqrt{Q}+(Q \ln Q)^{3/2}\right)$.

Аналогично, вычисляется и сумма $V_2$:
\begin{multline*}
V_2 = \sum_{\frac{Q}{4} < n\le \frac{Q}{2}} \frac{1}{4n} \sum_{|a|\le D} \sum_{-2n<c\le 2n} \sum_{2\sqrt{n}\le q\le Q} e^{2\pi i \frac{c(q^2-a)}{4n}} = \\
 =
\sum_{\frac{Q}{4} < n\le\frac{Q}{2}} \frac{1}{4n} \sum_{-2n<c\le 2n}
\left(\sum_{|a|\le D} e^{-2\pi i\frac{ac}{4n}}\right)\cdot \left(\sum_{2\sqrt{n}\le q\le Q} e^{2\pi i \frac{cq^2}{4n}}\right) =\\
=
\sum_{\frac{Q}{4}<n\le \frac{Q}{2}} \frac{1}{4n} \left(2[D]+1\right)\left([Q]-\left[2\sqrt{n}\right]+1\right) + \\
+
\sum_{\frac{Q}{4} < n\le\frac{Q}{2}} \frac{1}{4n} \sum_{\substack{-2n<c\le 2n\\c\ne 0}}
\left(\sum_{|a|\le D} e^{-2\pi i\frac{ac}{4n}}\right)\cdot \left(\sum_{2\sqrt{n}\le q\le Q} e^{2\pi i \frac{cq^2}{4n}}\right) = \\
=
\frac{DQ}{2} \sum_{\frac{Q}{4}<n \le \frac{Q}{2}} \frac1n + O\left(D\sqrt{Q} + (Q\ln Q)^{3/2}\right)
= \frac{\ln 2}{2}DQ + O\left(D\sqrt{Q} + (Q \ln Q)^{3/2}\right).
\end{multline*}

Окончательно находим:
\[
S_1'' = V_1 + V_2 = DQ\left(1+\frac{\ln 2}{2}\right) + O\left(D\sqrt{Q} + (Q \ln Q)^{3/2} \right).
\]

Далее,
\begin{multline*}
S_2'' = \sum_{1\le n\le \frac{Q}{2}} \sum_{|a|\le D} \sum_{\substack{2\sqrt{n} \le q \le \sqrt{2nQ}\\q^2\equiv a\!\!\!\!\!\pmod{4n}}} 1 = \\
=
\sum_{1\le n\le \frac{Q}{2}} \sum_{|a|\le D} \frac{1}{4n} \sum_{-2n<c\le 2n} \sum_{2\sqrt{n}\le q\le \sqrt{2nQ}} e^{2\pi i \frac{c(q^2-a)}{4n}} = \\
=
\sum_{1\le n\le \frac{Q}{2}} \frac{1}{4n} \sum_{-2n<c\le 2n}
\left(\sum_{|a|\le D} e^{-2\pi i \frac{ac}{4n}}\right)\cdot
\left(\sum_{2\sqrt{n}\le q\le \sqrt{2nQ}} e^{2\pi i\frac{cq^2}{4n}}\right) = \\
=
\sum_{1\le n\le \frac{Q}{2}} \frac{1}{4n} \left(2D\sqrt{2Qn} + O\left(\sqrt{Qn}+D\sqrt{n}\right)\right) + O\left((Q \ln Q)^{3/2}\right) = \\
=
\frac{D\sqrt{Q}}{\sqrt{2}} \sum_{1\le n\le\frac{Q}{2}}\frac{1}{\sqrt{n}} +
O\left(D\sqrt{Q} + (Q\ln Q)^{3/2}\right) =
DQ + O\left(D\sqrt{Q}+(Q \ln Q)^{3/2}\right).
\end{multline*}

Собирая вместе полученные оценки, найдём:
\begin{multline*}
N_1 = 2 S_1'' - S_2'' + O\left(D^{3/2}\ln Q + D^2/Q + D\sqrt{Q} + Q\right)= \\
=  DQ (1+\ln 2) + O\left(D^{3/2}\ln Q + D^2/Q+D\sqrt{Q}+(Q \ln Q)^{3/2}\right),
\end{multline*}
\[
N(Q,D) = 4 N_1 + O\left(Q\sqrt{D}+D^{3/2}\ln Q\right) = 4(1+\ln 2) DQ + O(r),
\]
где
\[
r = D^2/Q+D\sqrt{Q}+(Q \ln Q)^{3/2} + Q\sqrt{D}+D^{3/2}\ln Q = \sum_{k=1}^5 \Delta_k.
\]
Смысл обозначений $\Delta_k$, $1\le k\le 5$, очевиден. Замечая, что при $D\le 5Q^2$ имеет место неравенство
$\Delta_1 \ll \Delta_5$, при $D\le Q$ справедливы неравенство $\max(\Delta_2,\Delta_4)\ll \Delta_3$, а при $D > Q$ ---
неравенство $\max(\Delta_2,\Delta_4)\ll \Delta_5$, находим, что $r\ll \Delta_3 + \Delta_5$.

Формула для $N(Q,D)$ будет асимптотической, если $(Q \ln Q)^{3/2} + D^{3/2}\ln Q \ll DQ$, т.е.
$\sqrt{Q} (\ln Q)^{3/2} \ll D \ll \left(\frac{Q}{\ln Q}\right)^2$; изначально мы ограничивались значениями $D\le Q^2/2$.
Если положить $D= 5Q^{2-2v}$, то промежутку $\sqrt{Q} (\ln Q)^{3/2} \ll D \ll \left(\frac{Q}{\ln Q}\right)^2$ будет отвечать
промежуток $0 < v < \frac34$.

Теорема доказана.

\bigskip

{\small Ф.~Гётце}\\
{\footnotesize
{University of Bielefeld, 33501, Bielefeld, Germany}\\
E-mail: goetze@math.uni-bielefeld.de
}

\bigskip

{\small Д.\,В.~Коледа}\\
{\footnotesize
{Институт математики НАН Беларуси, Минск, Беларусь}\\
E-mail: koledad@rambler.ru
}

\bigskip

{\small М.\,А.~Королёв}\\
{\footnotesize
{МИАН им. В.А.~Стеклова, Москва, Россия}\\
E-mail: hardy\_ramanujan@mail.ru, korolevma@mi.ras.ru
}

\end{document}